# Stamps and Mathematics

**Nataliya M. Ivanova**

*Institute of Mathematics of NAS of Ukraine,*

*ivanova.nataliya@gmail.com*

***Abstract.*** *This study examines the potential of using math-themed postage stamps in mathematics lessons as a tool to engage students and integrate the subject with history, art, and culture. Since the first mathematical stamps appeared in the early 20th century, featuring prominent scholars like Carl Friedrich Gauss and Isaac Newton, they serve not only as philatelic artifacts but also as historical carriers of knowledge. The paper presents several practical projects to interest students, such as creating their own math stamps, investigating the price trends of math-themed stamps, and developing a timeline of mathematical discoveries depicted in philatelic issues. The proposed projects develop students' mathematical skills in areas such as percentage calculations, general arithmetic, working with time intervals, and statistical analysis. Students can analyze shapes, symmetry, and patterns on stamps, study principles of proportion, and explore geometric figures. Using stamps broadens students' horizons, providing an opportunity to become familiar with renowned mathematicians from different eras, countries, and cultures. This also offers students a new perspective on the subject, presenting mathematical discoveries as part of the world's cultural heritage. Postage stamps dedicated to mathematics can become a powerful tool for visualizing theoretical knowledge, stimulating interest in mathematics, and encouraging independent research among students[1]*

## 1. Introduction

Although philately has somewhat lost its popularity in recent years, it remains a fascinating activity that combines history, art, and culture. Collecting stamps offers the opportunity to explore global heritage, learn about notable personalities and events, and discover unique works of art. Moreover, philately fosters attention to detail, analytical thinking, and research skills. Reviving interest in philately could become an engaging and rewarding hobby for the modern generation, especially in the digital age, where the value of tangible, unique items is particularly cherished.

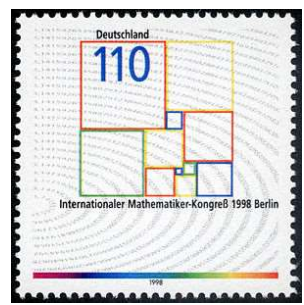

Philately, particularly collecting stamps with mathematical themes, is an engaging and effective tool for mathematics lessons. Mathematical stamps provide an opportunity to integrate mathematics with culture and history, helping students perceive it not merely as a set of formulas but as an integral part of the development of civilization. By examining stamps dedicated to prominent mathematicians, significant theorems, or groundbreaking discoveries,

---

students can learn more about the historical context in which these ideas were developed, thereby gaining a deeper appreciation for the importance of mathematical achievements.

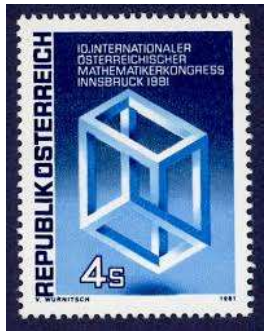 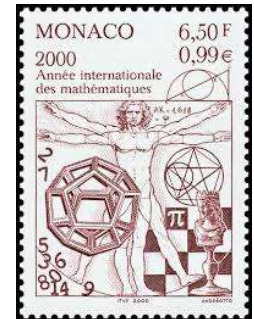

This approach can make mathematics more relatable and engaging for students, especially those with an inclination toward the humanities or arts. Using stamps as educational material develops visual perception and helps reveal the aesthetic side of mathematics through the exploration of symmetry, proportions, sequences, and other mathematical concepts. Moreover, working with or collecting stamps can spark students' research interest, inspire them to create their own projects, or even encourage the development of a new hobby.

With such a comprehensive approach, it becomes possible not only to enhance students' interest in mathematics but also to help them understand how mathematical ideas influence the world around them, even in such unexpected areas as philately.

## 2. History of mathematical philately

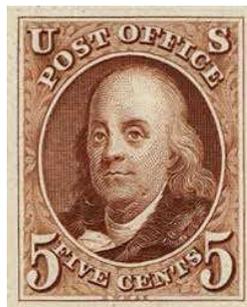

The history of mathematical stamps is a remarkable example of how mathematics can be represented in diverse forms, including philately. It is difficult to pinpoint exactly when this history began, as there is no universally accepted definition of what constitutes a mathematical stamp. Some collectors consider the 1847 United States issue depicting Benjamin Franklin to be the first mathematical stamp, even though Franklin was not primarily a mathematician (although his work in 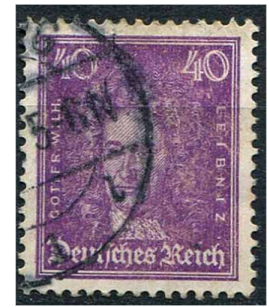 mathematical analysis and probability influenced the development of mathematics). Others point to the 1926 German stamp featuring Leibniz, whose primary profession as a mathematician is undisputed by either mathematicians or philatelists.

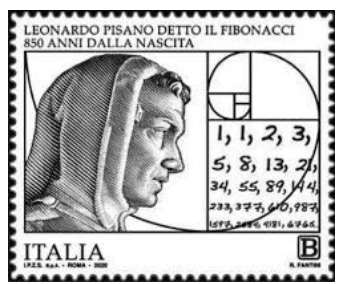

Mathematical stamps are postal issues that have a distinct connection to mathematics. They may feature renowned mathematicians or scientists who have made significant contributions to the development of mathematical science. These stamps often depict mathematical concepts or objects, formulas, symbols, or notations that have played a crucial role in this field of knowledge. Additionally, mathematical stamps frequently commemorate important events, such as national or international mathematical congresses or the International Mathematical Olympiad. Occasionally, they portray buildings dedicated to mathematical research or education, as well as the process of teaching mathematics itself.

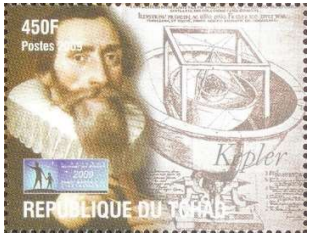
A separate category consists of stamps illustrating applied mathematics. These may relate to astronomy, ballistics, cartography, crystallography, navigation, mathematical art, games, or space exploration. Such stamps often depict objects or instruments connected to mathematics, such as astrolabes or slide rules. Their inclusion in mathematical collections depends on the collector's breadth of perspective.

However, not all materials can be classified as mathematical stamps. For example, this includes postal stationery, such as postcards featuring an image of a mathematician with a non-mathematical stamp affixed. Also excluded are so-called "Cinderella stamps", local or private issues that are not official. Non-postal stamps, such as revenue stamps, or stamps issued by non-existent or unrecognized countries are similarly excluded, particularly when their production exceeds actual postal needs.

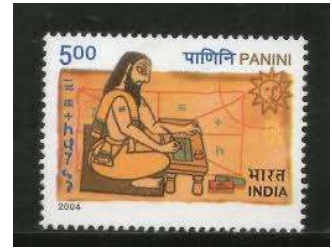

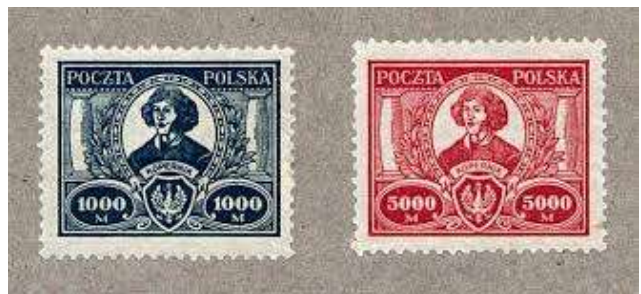

Since the 1920s, mathematical stamps have appeared regularly in various countries. Initially, they were dedicated to prominent mathematicians or significant scientific achievements closely linked to mathematics. The first depiction of a mathematician on postage stamps is considered to be the series issued to commemorate the 450th anniversary of Nicolaus Copernicus (Poland, 1923). Another highly popular series was the Norwegian stamps dedicated to Niels Abel. In 1955, Germany issued a stamp honoring Carl Friedrich Gauss, one of the most renowned mathematicians of all time, which significantly increased collectors' interest in mathematical stamps.

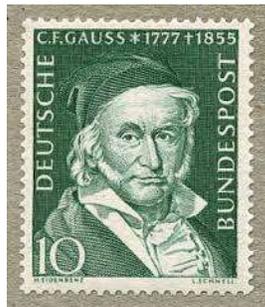

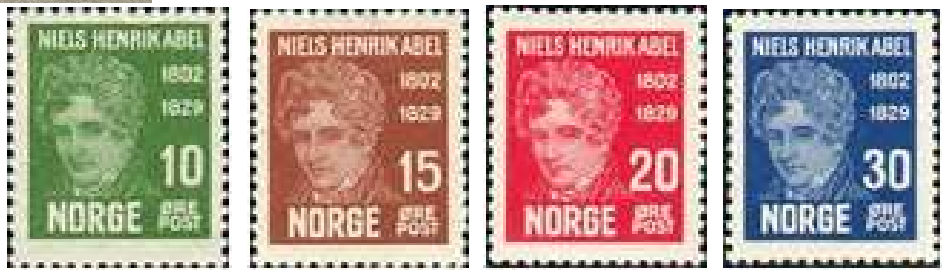

In the 1950s, countries began issuing stamps dedicated not only to individuals but also to significant mathematical ideas and theorems. One of the most famous examples is the Greek series from 1955, which honored Pythagoras and his theorem. In 1983, a stamp was issued to commemorate Euler's theorem, marking an important step in popularizing key mathematical concepts through philately. Fermat's Last Theorem is also frequently featured on stamps.

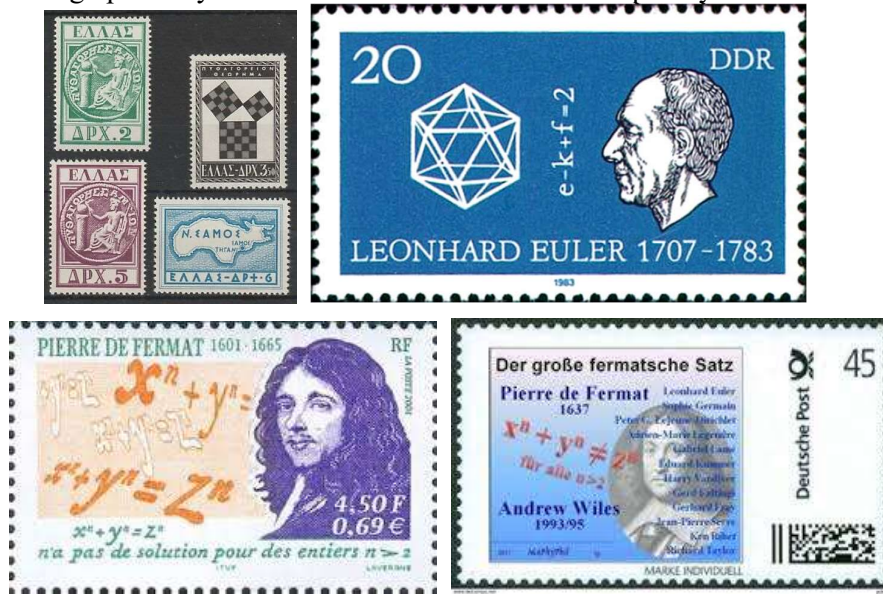

An important milestone was the series issued in 1971 in Nicaragua, which featured "The 10 Most Famous Mathematical Formulas."

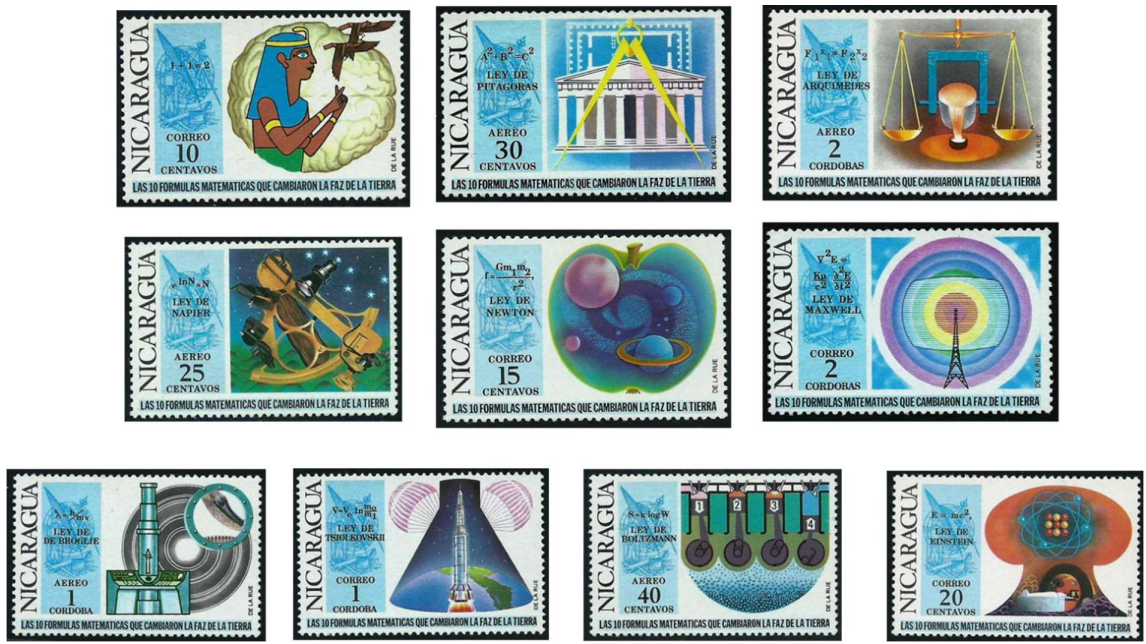

Of course, it is impossible to list all mathematical-themed stamps in a short article. Those interested can explore an extensive collection on the page [*Jeff Miller's postage stamps - MacTutor History of Mathematics*].

An exceptionally intriguing branch of mathematical philately is the collection of postmarks with mathematical themes.

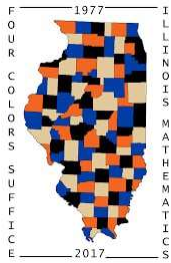
For a long time (1963–1971), the largest known Mersenne prime was the number $2^{11213} - 1$, which was discovered by Donald B. Gillies on a computer at the University of Illinois. To commemorate this achievement, the mathematics department of the university (which had its own post office) used a special postal cancellation for all its correspondence. This continued until around 1976, when Appel and Haken proved the Four Color Theorem, after which a new postal cancellation was created.

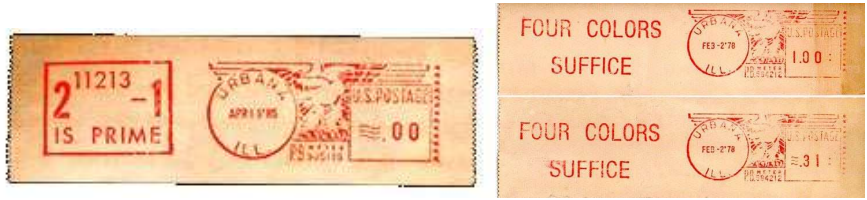

*(photo from open sources, author Chris Caldwell)*

In 2017, a commemorative, colourful postage mark was issued.

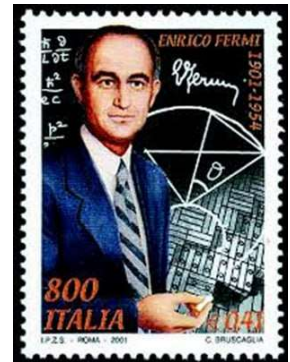
In 2017, a commemorative, colorful postage stamp was issued. Occasionally, interesting mistakes appear on stamps. For example, in 2001, to mark the 100th anniversary of Enrico Fermi's birth, the United States Postal Service decided to issue a stamp featuring his photograph. The scientist, who was known as "Papa" during his lifetime due to his almost flawless precision, was immortalized on a stamp with a 1948 photo. The irony is that there is a significant mistake in the photograph: the symbols "h" and "e" are swapped. The equation $a=h^2/(ec)$ should actually read $a=e^2/(hc)$. It's hard to believe that Fermi could have made such a mistake. Most likely not. It is probable that someone else wrote this equation on the board.

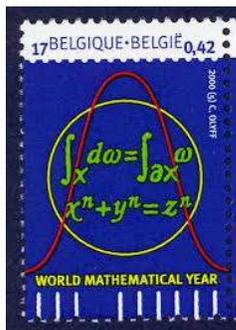
Mathematical stamps and postmarks have helped draw public attention to key mathematical ideas. Although their issuance has been limited, they have become an important element of philatelic collections worldwide. They can be viewed not only as monuments to major scientific achievements but also as tools for popularizing mathematical concepts among a broad audience. As a result, studying the history of mathematical stamps in mathematics lessons can help students not only learn more about the historical context of mathematical discoveries but also gain a better understanding of how mathematics influences culture and art. Collecting and studying such stamps can be an excellent way for students to see mathematics in a new light—as a part of cultural heritage that shapes our understanding of the world. Therefore, integrating mathematical stamps into the educational process not only enriches students' knowledge but also gives them the opportunity to see mathematical ideas as an important element of global scientific and cultural tradition.

### 3. Project activities of students

The history of mathematical stamps offers numerous opportunities for interesting and educational projects in mathematics lessons. By studying these stamps, students can not only learn about prominent mathematicians and their achievements but also apply their mathematical knowledge in real-world projects that help develop their skills. Examples of such projects are provided in Table 1.

Table 1

Examples of students' research projects

| Project | Short description | Mathematical skills |
|---|---|---|
| Creating a Mathematical Stamp Collection | Students can research stamps from different countries dedicated to prominent mathematicians, theorems, or mathematical models. Then, they can create a presentation explaining the significance of each stamp and, for example, calculate how many math-themed stamps each country has issued. | counting, percentage calculations, comparisons, compiling statistics |
| Mathematicians' Birthdays | Distribute images of mathematicians featured on stamps among the students and ask them to find the birthdate of each mathematician. Based on this information, students can create a timeline with mathematicians from different eras. | working with time intervals, chronology, creating charts and graphs |
| Creating a Personal Math-Themed Stamp | Ask students to design a stamp with a mathematical theme: famous formulas, figures, models, or mathematical symbols. Additionally, they can be asked to come up with a nominal value for the stamp and calculate how its value would change over time, considering inflation. | working with percentages, basic economics, creative skills |
| Mathematical Shapes and Symmetry on Stamps | Students will explore the geometric shapes and symmetry used in stamp designs. They can analyse proportions, symmetry of images, repeating elements, and even identify patterns related to Fibonacci numbers or other mathematical sequences. | geometry, symmetry, sequences, patterns |
| Calculating Value and Investments | Researching the historical changes in the value of famous stamps dedicated to mathematicians or mathematical themes will help students understand how collectible value increases. Students can try to predict the price of a stamp in a few years. | percentages, statistics, creating graphs, working with time |
| Mathematical Quest "Find Mathematics on Stamps" | Prepare a collection of various stamps with hidden mathematical symbols, shapes, or formulas. Students need to find all the mathematical elements on each stamp and | pattern recognition, analytical thinking, knowledge of |

| | explain how they are related to mathematical concepts. | different areas of mathematics |
|---|---|---|
| Creating a "Mathematical Postal History" | Students can create a chronology of mathematical discoveries, each represented by a stamp or postcard depicting a particular era in mathematics. This will help them study the development of mathematical ideas from ancient times to the present day. | working with numbers and time, the history of mathematical concepts, creating timelines |

Such projects stimulate students' interest in the subject, while also developing critical thinking, research skills, and the ability to work with various types of information, which are essential for forming a broad worldview.